\documentclass{article}
\usepackage{amsmath}
\usepackage{amssymb}
\usepackage{amsfonts}
\usepackage{amsthm}

\DeclareMathOperator{\vol}{vol}
\DeclareMathOperator{\lin}{lin}
\DeclareMathOperator{\interior}{int}
\DeclareMathOperator{\conv}{conv}

\title{A Discrete Analogue for Minkowski's Second Theorem on Successive Minima}
\author{Romanos-Diogenes Malikiosis}
\date{July 5th, 2009}

\begin{document}
\maketitle

\begin{abstract}
\noindent
The main result of this paper is an inequality relating the lattice point enumerator of a 3-dimensional, 0-symmetric convex body and its successive minima. This is an example of generalization of Minkowski's theorems on successive minima, where the volume is replaced by the discrete analogue, the lattice point enumerator. This problem is still open in higher dimensions, however, we introduce a stronger conjecture that shows a possibility of proof by induction on the dimension.
\end{abstract}

\section{Introduction}

\noindent
A subset $K$ of $\mathbb{R}^d$ is called a convex body if it is convex, compact, with nonempty interior (sometimes the technical condition $K\subset \overline{\interior(K)}$ is required, but we will not need it here). The set of all convex bodies of $\mathbb{R}^d$ will be denoted as $\mathcal{K}^d$, and the subset of the $0$-symmetric elements will be denoted as $\mathcal{K}_0^d$. Furthermore, $\mathcal{L}^d$ will denote the set of all lattices in $\mathbb{R}^d$. We will denote the $i$th successive minimum of $K\in \mathcal{K}_0^d$ with respect to $\Lambda\in\mathcal{L}^d$ by $\lambda_i(K,\Lambda)$, i.e.,
\[\lambda_i=\inf\{\lambda>0|\dim(\lambda K\cap\Lambda) \geqslant i\},\]
where $\dim(A)$ for $A\subset\mathbb{R}^d$ is the dimension of the vector space spanned by all vectors of the set $A$.

\noindent
Minkowski proved the following two inequalities relating the volume of $K\in\mathcal{K}_0^d$, i.e., its $d$-dimensional Lebesgue measure, with its successive minima, with respect to a lattice $\Lambda$:
\begin{equation}
\frac{\vol(K)}{\det(\Lambda)}\leqslant \left(\frac{2}{\lambda_1(K,\Lambda)}\right)^d
\end{equation}
and
\begin{equation}
\frac{1}{d!}\prod_{i=1}^d \frac{2}{\lambda_i(K,\Lambda)}\leqslant \frac{\vol(K)}{\det(\Lambda)}\leqslant\prod_{i=1}^d \frac{2}{\lambda_i(K,\Lambda)}
\end{equation}
which is an improvement of (1). The above are known as the first and second theorem on successive minima, respectively. About a century later, in 1993, Betke, Henk, and Wills [BHW] stated discrete analogues of these theorems, where the volume is replaced by the lattice point point enumerator, $G(K,\Lambda):=\#(K\cap\Lambda)$. They proved the analogue for the first theorem, which predicts the following inequality:
\begin{equation}
G(K,\Lambda)\leqslant \left[\frac{2}{\lambda_1(K,\Lambda)}+1\right]^d
\end{equation}
As for the second theorem, they stated a conjecture, which they verified in the planar case. The conjectural inequality is:
\begin{equation}
G(K,\Lambda)\leqslant \prod_{i=1}^d \left[\frac{2}{\lambda_i(K,\Lambda)}+1\right]
\end{equation}
In section 3.2 we shall see a proof for the $3$-dimensional case of the above conjecture. There is a notion of induction in the proof; we need statements about intersections of a given convex body by hyperplanes passing through lattice points. The resulting convex bodies, whose dimension is $d-1$, are not centrally symmetric, in general. Therefore, it is necessary to extend the definition of the successive minima, as well as the results referred to them, to the class of all convex bodies, not necessarily $0$-symmetric, namely
\[\lambda_i(K,\Lambda):=\lambda_i(\tfrac{1}{2}\mathfrak{D}K,\Lambda),\]
where $\mathfrak{D}K$ stands for the difference body of $K$, i.e.,
\[\mathfrak{D}K=K-K=\{x-y|x,y\in K\}.\]
Under this definition, Minkowski's theorems still hold; this is a simple consequence of the Brunn-Minkowski theorem ([GL],pp. 12 and 32), which predicts that
\[\vol(K)\leqslant \vol(\tfrac{1}{2}\mathfrak{D}K).\]
In section 3.1, we will provide a proof for inequality (3) for the non $0$-symmetric case; we conjecture that inequality (4) still holds when $K$ is not $0$-symmetric, and we will prove some weaker estimates, exactly as in [M].

\section{Some auxiliary lemmata}

\noindent
We will first need the following standard convention; remember that the successive minima of $K\in\mathcal{K}^d$ with respect to a lattice $\Lambda$ are those of $\frac{1}{2}\mathfrak{D}K$. By definition of the successive minima $\lambda_i(K,\Lambda)$, there are $d$ linearly independent lattice vectors $a^i=a^i(K,\Lambda)$ such that
\[a^i(K,\Lambda)\in \tfrac{\lambda_i(K,\Lambda)}{2}\mathfrak{D}K\cap\Lambda.\]
Then, we construct a basis of $\Lambda$, denoted by $e^i=e^i(K,\Lambda)$, $1\leqslant i\leqslant d$, such that
\[\lin(a^1,\dotsc,a^i)=\lin(e^1,\dotsc,e^i)\]
for all $i$, $1\leqslant i\leqslant d$. Furthermore, we define the following subgroups of $\Lambda$:
\[\Lambda^i(K):=\mathbb{Z}e^1\oplus\dotsb\oplus\mathbb{Z}e^i\]
We will usually abbreviate the notation to $\Lambda^i$. It should be noted that there is an abuse of notation here; it is evident that the choice of the $a^i$'s and the $e^i$'s, as well as the $\Lambda^i$'s, is not always unique. However, by this notation we shall always mean a choice of vectors or subgroups with the above properties. The main property that will be used later is
\begin{equation}
\interior(\tfrac{\lambda_i}{2}\mathfrak{D}K)\cap\Lambda\subset\Lambda^{i-1}
\end{equation}

\noindent
\textbf{\underline{Lemma 2.1}} \textit{Let $K\in\mathcal{K}^d$, $\Lambda\in\mathcal{L}^d$. For each real $n_i$, satisfying $n_i>2/\lambda_i$, we have
\[\mathfrak{D}K\cap n_i(\Lambda\setminus \Lambda^{i-1})=\varnothing\].}

\noindent
\textbf{\underline{Proof}:} Assume otherwise; then the intersection
\[\tfrac{1}{n_i}\mathfrak{D}K\cap(\Lambda\setminus\Lambda^{i-1})\]
would be nonempty. The left part of this intersection is a subset of
\[\interior(\tfrac{\lambda_i}{2}\mathfrak{D}K),\]
since $n_i>2/\lambda_i$. Therefore, the intersection
\[\interior(\tfrac{\lambda_i}{2}\mathfrak{D}K)\cap(\Lambda\setminus\Lambda^{i-1})\]
would be nonempty; but this is absurd, since it contradicts (5) above.$\square$

\noindent
The following is an adaptation of lemma 2.1 in [H], for the case of non-symmetric convex bodies. The proof is identical, so we will not provide it here.

\noindent
\textbf{\underline{Lemma 2.2}} \textit{Let $K\in\mathcal{K}^d$, $\Lambda\in\mathcal{L}^d$, and $\widetilde{\Lambda}$ a sublattice of $\Lambda$. Then
\[G(K,\Lambda)\leqslant \frac{\det(\widetilde{\Lambda})}{\det(\Lambda)}G(\mathfrak{D}K,\widetilde{\Lambda}).\]
}

\noindent
The following two lemmata will be used for the proof of inequality (4) in the $3$-dimensional case. Notice that they are statements in $d$ dimensions.

\noindent
\textbf{\underline{Lemma 2.3}} \textit{Let $K\subset\mathbb{R}^d$ be a convex body, $\Lambda\in\mathcal{L}^d$, such that $K\cap\Lambda=\varnothing$. For a real $t>1$, there is some $v\in\Lambda$ such that
\[K\cap(v+t\Lambda)=\varnothing.\]
}

\noindent
\textbf{\underline{Proof}:} Take $v\in\Lambda$ such that $\#(\conv(v,K)\cap\Lambda)$ is minimal. If this number is greater than 1, then there is some $w\in\Lambda$, $w\neq v$, such that $w\in\conv(v,K)$. Hence, $\conv(w,K)\subset\conv(v,K)$, and $v\notin\conv(w,K)$, contradicting the minimality of $\#(\conv(v,K)\cap\Lambda)$. Thus, $\conv(v,K)\cap\Lambda=\{v\}$. We claim that $K\cap(v+t\Lambda)=\varnothing$. Assume not; then there is some $u\in\Lambda$ such that $v+tu\in K$. By convexity, and the fact that $t>1$, we get $v+u\in\conv(v,K)$, which implies $u=0$, or $v\in K$, a contradiction, since $K\cap\Lambda=\varnothing$. This concludes the proof.$\square$

\noindent
The next lemma generalizes the above:

\noindent
\textbf{\underline{Lemma 2.4}} \textit{Let $K\subset\mathbb{R}^d$ be a convex body, and $\Lambda\in\mathcal{L}^d$. Let $S\subset\Lambda$ be finite, and $r$ be a positive integer, such that
\begin{description}
	\item[(1)] $(K-S)\cap r\Lambda=\varnothing$.
	\item[(2)] $\mathfrak{D}S\cap r(\Lambda\setminus\{0\})=\varnothing$. 
\end{description}
Now let $t>r$ be an integer. There is a set $S'\subset\Lambda$, obtained by translating each $v\in S$ by some vector $r\cdot w(v)$, where $w(v)\in\Lambda$, such that
\begin{description}
	\item[(1)$'$] $(K-S')\cap t\Lambda=\varnothing$.
	\item[(2)$'$] $\mathfrak{D}S'\cap t(\Lambda\setminus\{0\})=\varnothing$. 
\end{description}
}

\noindent
\textbf{\underline{Proof}:} We will prove it by induction on $\#(S)$. If $\#(S)=1$, i.e. $S=\{v\}$, we use Lemma 2.3 for $K-v$ and the lattice $r\Lambda$. Since $t>r$, there is some $w(v)\in\Lambda$, such that $(K-v)\cap(r\cdot w(v)+t\Lambda)=\varnothing$. Put $S'=\{v+r\cdot w(v)\}$, and we can easily see that $(1)'$ is satisfied. It should be noted that when $\#(S)=1$, conditions $(2)$ and $(2)'$ hold vacuously.

\noindent
Now, assume that $\#(S)>1$. Take $v\in S+r\Lambda$, such that $\#(\conv(v,K)\cap(S+r\Lambda))$ is minimal. Again, as in the proof of Lemma 2.3, we must have $\conv(v,K)\cap(S+r\Lambda)=\{v\}$. Apply induction for $\widetilde{K}=\conv(v,K)$ and $\widetilde{S}=S\setminus(S\cap(v+r\Lambda))$; we have $\#(\widetilde{S})=\#(S)-1$. Let's see why (1) and (2) are satisfied for $\widetilde{K}, \widetilde{S}$ (same $r,\Lambda$); (2) is obviously satisfied, as $\widetilde{S}\subset S$. If (1) were not satisfied, there would be some $w\in\widetilde{S}$ and $u\in\Lambda$ such that $w+ru\in\conv(v,K)$. By minimality assumption, $w+ru=v$. But $v\notin \widetilde{S}+r\Lambda$, a contradiction. Thus, (1) and (2) hold for $\widetilde{K}$, $\widetilde{S}$, and by induction there is some $\widetilde{S}'\subset\Lambda$, obtained from $\widetilde{S}$ by translating each $u\in\widetilde{S}$ by $r\cdot w(u)$, $w(u)\in\Lambda$, such that
\[(\widetilde{K}-\widetilde{S}')\cap t\Lambda=\varnothing,\]
and
\[\mathfrak{D}\widetilde{S}'\cap t(\Lambda\setminus\{0\})=\varnothing.\]
Now, simply take $S'=\widetilde{S}'\cup\{v\}$. (2)$'$ is satisfied for $S'$; if $x,y\in\widetilde{S}'$, then $x-y\notin t(\Lambda\setminus\{0\})$ from the above. If $x\in\widetilde{S}'$ and $y=v$, then again from the above, $v-x\notin t\Lambda$, since $v\in K'$. If $x=y=v$, we have nothing to prove, so
\[\mathfrak{D}S'\cap t(\Lambda\setminus\{0\})=\varnothing.\]
(1)$'$ is also satisfied for $K,S'$; assume not. Then, there would be some $w\in S'$, $u\in\Lambda$ such that $w+tu\in K$. If $w\in\widetilde{S}'$, then $w+tu\in \widetilde{K}$, which contradicts $(\widetilde{K}-\widetilde{S}')\cap t\Lambda=\varnothing$. If $w=v$, then $v+tu\in K$, and by convexity, $v+ru\in K$, hence $u=0$, by minimality assumption, and $v\in K$, a contradiction. This concludes the proof.$\square$

\section{Inequalities for $G(K,\Lambda)$}

\noindent
Throughout the rest of the paper we will use the notation
\[q_i(K,\Lambda)=\left[\frac{2}{\lambda_i(K,\Lambda)}+1\right].\]
Also, when $\Lambda$ is the standard lattice $\mathbb{Z}^d$, we write $G(K)$ instead of $G(K,\mathbb{Z}^d)$.

\subsection{The general case}

\noindent
The method of the following proof is similar to the proof of theorem 1.5 in [H], from a slightly different viewpoint.

\noindent
\textbf{\underline{Theorem 3.1.1}} \textit{Let $K\in\mathcal{K}^d$, $\Lambda\in\mathcal{L}^d$, $q_i=q_i(K,\Lambda)$. Let also $n_1,\dotsc,n_d$ be a sequence of integers satisfying
\begin{itemize}
	\item $n_{i+1}$ divides $n_i$, $1\leqslant i\leqslant d-1$.
	\item $q_i\leqslant n_i$, $1\leqslant i \leqslant d$.
\end{itemize}
Then,
\[G(K,\Lambda)\leqslant \prod_{i=1}^{d}n_i.\]}

\noindent
\textbf{\underline{Proof}:} Let $e^i=e^i(K,\Lambda)$ and define
\[\widetilde{\Lambda}=\mathbb{Z}n_1e^1\oplus\dotsb\oplus\mathbb{Z}n_de^d.\]
By lemma 2.2,
\[G(K,\Lambda)\leqslant \frac{\det(\widetilde{\Lambda})}{\det(\Lambda)}G(\mathfrak{D}K,\widetilde{\Lambda})= G(\mathfrak{D}K,\widetilde{\Lambda})\prod_{i=1}^{d}n_i.\]
It suffices to prove that $G(\mathfrak{D}K,\widetilde{\Lambda})=1$, or equivalently
\[\mathfrak{D}K\cap(\widetilde{\Lambda}\setminus\{0\})=\varnothing.\]
This follows from lemma 2.1 and the fact that
\[\widetilde{\Lambda}\setminus\{0\}\subset \bigcup_{i=1}^{d}n_i(\Lambda\setminus\Lambda^{i-1})\]
(we remind that $n_i\geqslant q_i>2/\lambda_i$). Indeed, let $g\in\widetilde{\Lambda}\setminus\{0\}$ be arbitrary, and let $k$ be minimal such that
\[g\in\mathbb{Z}n_1e^1\oplus\dotsb\oplus\mathbb{Z}n_ke^k.\]
Since $n_k$ divides all $n_1,\dotsc,n_{k-1}$ by hypothesis, we have $g\in n_k\Lambda$. By minimality of $k$, we also have $g\notin \Lambda^{k-1}$, hence $g\in n_k(\Lambda\setminus\Lambda^{k-1})$, as desired.$\square$

\noindent
We obtain inequality (3) for the non $0$-symmetric case as well by simply putting $n_i=q_1$, for all $i$, $1\leqslant i\leqslant d$. We remind the following definition, given in [M]:

\noindent
\textbf{\underline{Definition 3.1.1}} \textit{Let $C_d$ denote the least positive constant, such that for any sequence of $d$ positive integers, $x_1\leqslant x_2,\leqslant\dotsb\leqslant x_d$, there exists a sequence of integers $y_1,y_2,\dotsc,y_d$ satisfying:
\begin{description}
	\item[a.] $x_i\leqslant y_i$, for all $i$, $1\leqslant i\leqslant d$.
	\item[b.] $y_i$ divides $y_{i+1}$, for all $i$, $1\leqslant i\leqslant d-1$.
	\item[c.] $\displaystyle\frac{y_1y_2\dotsm y_d}{x_1x_2\dotsm x_d}\leqslant C_d$.
\end{description}}

\noindent
We can always choose $n_1,\dotsc,n_d$ with the properties given in theorem 3.1.1, such that
\[\prod_{i=1}^{d}n_i\leqslant C_d\prod_{i=1}^{d}q_i.\]
The estimates for $C_d$ in [M], yield the following:

\noindent
\textbf{\underline{Theorem 3.1.2}} \textit{Let $K\in\mathcal{K}^d$, $\Lambda\in\mathcal{L}^d$. Then
\[G(K,\Lambda)\leqslant \frac{4}{e}(\sqrt{3})^{d-1}\prod_{i=1}^{d}q_i(K,\Lambda).\]
If $K\in\mathcal{K}_0^d$, then
\[G(K,\Lambda)\leqslant \frac{4}{e}\left(\sqrt[3]{\frac{40}{9}}\right)^{d-1}\prod_{i=1}^{d}q_i(K,\Lambda).\]}

\noindent
Notice that when $K$ is non 0-symmetric, we cannot disregard all $q_i$ that are less than or equal to $2$, as it was done in [M]; the reason is that we may have $q_1=\dotsb=q_d=2$, and $K\cap\Lambda$ have full affine dimension (while if $K$ was 0-symmetric, this would mean that $K\cap\Lambda=\{0\}$). Consider, for example, the $d$-dimensional cube $[0,1]^d$. So, in the non 0-symmetric case we cannot reduce the base to $(\sqrt[3]{40/9})^{d-1}$, using this argument.

\subsection{The case $d=3$}

\noindent
We will introduce an inductive method in order to prove inequality (4); this method works up to dimension 3, and it is inadequate for higher dimensions. As we will see in the next section, stronger versions of lemmata 2.3 and 2.4 might be needed for this method to work in all dimensions.

\noindent
Let $K\in\mathcal{K}_0^d$, $\Lambda\in\mathcal{L}^d$. Fix a basis $e^i=e^i(K,\Lambda)$ of $\Lambda$, that satisfies the properties given in section 2. We will write each vector $x$ of $\mathbb{R}^d$ with coordinates with respect to this basis:
\begin{eqnarray*}
x &=& (x_1,\dotsc,x_d)\\
&=& x_1e^1+\dotsb+x_de^d.
\end{eqnarray*}
Define
\[K[t]:=\{x\in K|x_d=t\},\]
i.e., the subset of $K$ whose elements have fixed height, or the intersection of $K$ by a hyperplane parallel to the vector subspace spanned by $e^1,\dotsc,e^{d-1}$. We can write $G(K,\Lambda)$ in terms of lattice point enumerators of convex bodies whose dimension is $d-1$; this is the point where induction could be used. Namely,
\[G(K,\Lambda)=\sum_{t\in\mathbb{Z}}G(K[t]-te^d,\Lambda^{d-1}).\]
The bodies $K[t]-te^d$ are projections of the intersections $K[t]$ on the vector subspace spanned by $e^1,\dotsc,e^{d-1}$ along the lattice vector $e^d$. As before, $\Lambda^{d-1}$ is the $\mathbb{Z}$-span of $e^1,\dotsc,e^{d-1}$. Apart from $K[0]$ which is $0$-symmetric, the other projections are not necessarily $0$-symmetric. This is the main reason for extending inequalities (3) and (4) for the non symmetric case.

\noindent
Next, we observe that
\[\tfrac{1}{2}\mathfrak{D}(K[t]-te^d)\subset\tfrac{1}{2}\mathfrak{D}K,\]
therefore, for $1\leqslant i\leqslant d-1$
\[\lambda_i(K[t]-te^d,\Lambda^{d-1})\geqslant\lambda_i(K,\Lambda),\]
which implies
\[q_i(K[t]-te^d,\Lambda^{d-1})\leqslant q_i(K,\Lambda),\]
for $1\leqslant i\leqslant d-1$. Assuming that inequality (4) holds for $d-1$, we would have
\[G(K[t]-te^d,\Lambda^{d-1})\leqslant\prod_{i=1}^{d-1}q_i(K,\Lambda),\]
for all $t\in\mathbb{Z}$. Only the factor $q_d$ is missing; normally, it would be obtained from the number of the nonempty intersections, $K[t]$. But it is not always the case that this number is less than $q_d(K,\Lambda)$.

\noindent
The next technique is to group all intersections whose heights are congruent modulo $q_d$. Then, the above sum becomes
\[G(K,\Lambda)=\sum_{r=0}^{q_d-1}\sum_{t\equiv r(\text{mod}q_d)}G(K[t]-te^d,\Lambda^{d-1}).\]
It suffices to prove that for each fixed $r$, we have
\[\sum_{t\equiv r(\text{mod}q_d)}G(K[t]-te^d,\Lambda^{d-1})\leqslant \prod_{i=1}^{d-1}q_i.\]
Of course, we could have more than one convex body in the above sum, however, the above collection of convex bodies $K[t]-te^d$, $t\equiv r(\text{mod}q_d)$ satisfies some restricting conditions, namely
\begin{description}
	\item[(1)] $\mathfrak{D}(K[t]-te^d)\cap q_i(\Lambda^{d-1}\setminus\Lambda^{i-1})=\varnothing$ for all $t\equiv r(\text{mod}q_d)$ and $1\leqslant i\leqslant d-1$.
	\item[(2)] $((K[t]-te^d)-(K[t']-t'e^d))\cap q_d\Lambda^{d-1} = \varnothing$ for all $t,t'\equiv r(\text{mod}q_d)$, $t\neq t'$. 
\end{description}
The above are consequences of lemma 2.1. Indeed, for the first one we observe that
\[\mathfrak{D}(K[t]-te^d)\cap q_i(\Lambda^{d-1}\setminus\Lambda^{i-1})\subset\mathfrak{D}K\cap q_i(\Lambda\setminus\Lambda^{i-1}),\]
and the latter is empty since $q_i>2/\lambda_i$. As for the second, if
\[((K[t]-te^d)-(K[t']-t'e^d))\cap q_d\Lambda^{d-1} \neq\varnothing,\]
there would exist some $v\in\Lambda^{d-1}$ such that $q_dv+(t-t')e^d\in K[t]-K[t']\subset \mathfrak{D}K$. Since $q_d|t-t'$, and $t\neq t'$, the intersection
\[\mathfrak{D}K\cap q_d(\Lambda\setminus\Lambda^{d-1})\]
would be nonempty, again a contradiction by lemma 2.1.

\noindent
It is natural to state the following conjecture:

\noindent
\textbf{\underline{Conjecture 3.2.1}} \textit{Let $K_1,\dotsc,K_n\subset \mathbb{R}^d$ be convex bodies and $\Lambda\in\mathcal{L}^d$. Also, let $e^1,\dotsc,e^d$ be a basis of $\Lambda$, and denote by $\Lambda^i$ the $\mathbb{Z}$-span of $0,e^1,\dotsc,e^i$, and let $q_1\geqslant q_2\geqslant \dotsb \geqslant q_{d+1}$ be positive integers satisfying
\begin{description}
	\item[(1)] $\mathfrak{D}K_j\cap q_i(\Lambda\setminus\Lambda^{i-1})=\varnothing$ for all $1\leqslant j\leqslant n$ and $1\leqslant i\leqslant d$.
	\item[(2)] $(K_j-K_l)\cap q_{d+1}\Lambda = \varnothing$ for all $1\leqslant j,l\leqslant n$, $j\neq l$. 
\end{description}
Then
\[\sum_{j=1}^n G(K_j,\Lambda)\leqslant \prod_{i=1}^d q_i.\]
}

\noindent
From the above analysis, it is clear that the above conjecture implies inequality (4) for one dimension higher. We will verify this conjecture for $d=1,2$, thus proving inequality (4) in all dimensions up to three. A statement in support of this conjecture is that condition (2) is too restricting for the convex bodies $K_j$, given the fact that $q_{d+1}$ is smaller than the rest of the $q_i$'s. It simply says that no two translates of $K_j$ and $K_l$, $j\neq l$, by vectors of $q_{d+1}\Lambda$ intersect. In the next section, we will see a reduction of this conjecture, which is more convincing.

\noindent
\textit{\underline{Proof of Conjecture 3.2.1, $d=1$}:} Without loss of generality, we may assume that $\Lambda=\mathbb{Z}$. Let $K_j=[a_j,b_j]$, $1\leqslant j\leqslant n$. Conditions (1) and (2) read
\begin{description}
	\item[(1)] $b_j-a_j<q_1$ for all $1\leqslant j\leqslant n$.
	\item[(2)] $(K_j-K_l)\cap q_2\mathbb{Z} = \varnothing$ for all $1\leqslant j,l\leqslant n$, $j\neq l$. 
\end{description}
If $b_1-a_1\geqslant q_2$, then the union of $K_1$ with all its translates by multiples of $q_2$ cover all of $\mathbb{R}$, so by condition (2) we must have $n=1$, therefore
\[\sum_{j=1}^{n}G(K_j,\Lambda)=G(K_1)\leqslant q_1\]
by (1). If $b_1-a_1<q_2$, there is a translate of each $K_j$ by some multiple of $q_2$, $2\leqslant j\leqslant n$, that lies in $(b_1,a_1+q_2)$, again by (2). Since they do not intersect each other by (2), we will have
\[\sum_{j=1}^{n}G(K_j)\leqslant G([a_1,a_1+q_2))=q_2\leqslant q_1.\square\]

\noindent
\textit{\underline{Proof of Conjecture 3.2.1, $d=2$}:} Let 
\[D=\dim\biggl(\biggl(\bigcup_{j=1}^n \mathfrak{D}K_j\biggr)\cap q_3\Lambda\biggr).\]
We distinguish cases for $D$:

\noindent
\underline{$D\leqslant 1$}: There exists a primitive lattice vector, say $v$, such that
\[\biggl(\bigcup_{j=1}^n \mathfrak{D}K_j\biggr)\cap q_3\Lambda\subset\mathbb{Z}(q_3v)\]
therefore
\[\bigl(\bigl(\bigcup_{j=1}^n\mathfrak{D}K_j\bigr)\cap q_3\bigl(\Lambda\setminus\mathbb{Z}v\bigr)\bigr)=\varnothing.\]
Find $w\in \Lambda$ such that $v,w$ is a basis for $\Lambda$. Then
\[\sum_{j=1}^nG(K_j,\Lambda)=\sum_{r=0}^{q_3-1}\sum_{j=1}^n\sum_{t\equiv r(\text{mod}q_3)}G(K_j[t]-tw,\mathbb{Z}v).\]
To prove that the above sum is less than or equal to $q_1q_3$ it suffices to prove that
\[\sum_{j=1}^n\sum_{t\equiv r(\text{mod}q_3)}G(K_j[t]-tw,\mathbb{Z}v)\leqslant q_1,\]
for a fixed $r$, where the notation $K_j[t]$ refers to the basis $v,w$. Naturally, we identify $\mathbb{R}v$ with $\mathbb{R}$, so the collection of all sets $K_{j,t}:=K_j[t]-tw$ (which only finitely of them are nonempty) is a collection of compact intervals on $\mathbb{R}$. We have
\[\mathfrak{D}K_{j,t}\cap q_1(\mathbb{Z}v\setminus\{0\})\subset \mathfrak{D}K_j\cap q_1(\Lambda\setminus\{0\})=\varnothing,\]
for all $j$, so condition (1) is satisfied. Furthermore, when $t\neq t'$, if the intersection
\[(K_{j,t}-K_{j,t'})\cap q_3(\mathbb{Z}v)\]
is nonempty, then there exists $u\in \mathbb{Z}v$ such that
\[q_3u+(t-t')w\in K_j[t]-K_j[t']\subset \mathfrak{D}K_j,\]
which implies that
\[\mathfrak{D}K_j\cap q_3(\Lambda\setminus\mathbb{Z}v)\neq \varnothing,\]
a contradiction, since $D\leqslant 1$. If $i\neq j$, and if the intersection
\[(K_{i,t}-K_{j,t'})\cap q_3(\mathbb{Z}v)\]
is nonempty, then there is $u\in \mathbb{Z}v$ such that
\[q_3u+(t-t')w\in K_i[t]-K_j[t']\subset K_i-K_j,\]
which implies that
\[(K_i-K_j)\cap q_3\Lambda\neq\varnothing,\]
again a contradiction. So, condition (2) is satisfied, and since the 1-dimensional case is true, we will have
\[\sum_{j=1}^n\sum_{t\equiv r(\text{mod}q_3)}G(K_j[t]-tw,\mathbb{Z}v)\leqslant q_1,\]
as desired.

\noindent
\underline{$D=2$}: This means that there are two primitive, linearly independent vectors of $\Lambda$ in $\bigcup\mathfrak{D}K_j$, say $v,w$. We may assume that $q_3v\in\mathfrak{D}K_i$ and $q_3w\in\mathfrak{D}K_j$, for some indices $i,j$. We want to show that $i=j$ (if $n=1$, this is vacuously true, so we assume $n\geqslant2$). We have
\[K_i\cap(K_i-q_3v)\neq\varnothing,\]
so we pick an element $x$ from this intersection. Hence, $x,x+q_3v\in K_i$. Also,
\[K_j\cap(K_j+q_3w)\neq\varnothing,\]
and we pick an element $y$ from this intersection, so $y,y-q_3w\in K_j$. Let $\widetilde{\Lambda}=\mathbb{Z}v\oplus\mathbb{Z}w$, and consider the fundamental parallelogram of $q_3\widetilde{\Lambda}$ with vertices $x,x+q_3v,x+q_3w,x+q_3(v+w)$, say $\mathcal{P}$. Since $\mathcal{P}$ is a fundamental parallelogram, there is a translate of $y$ by $q_3\widetilde{\Lambda}$ (and hence by $q_3\Lambda$ as well) in $\mathcal{P}$. Without loss of generality, we may assume that $y\in\mathcal{P}$ (if we translate any $K_i$ by an element of $q_3\Lambda$, conditions (1) and (2) still hold). So, assume that $y=x+\alpha q_3v+\beta q_3w$, where $0\leqslant \alpha,\beta<1$. But then the element
\[y-\beta q_3w=x+\alpha q_3v\]
belongs to both $\conv(x,x+q_3v)$ and $\conv(y,y-q_3w)$, i.e., the intersection $K_i\cap K_j$ is nonempty. This contradicts condition (2), if $i\neq j$, so we must have $i=j$.

\noindent
Without loss of generality, assume that $i=1$, that is, $v,w\in \mathfrak{D}K_1$. Choose $v,w$ so that the index $[\Lambda :\widetilde{\Lambda}]$ is minimal. If $[\Lambda :\widetilde{\Lambda}]>1$, then there is a point of $q_3\Lambda$ in $\conv(0,q_3v,q_3w)$, different from $0,q_3v,q_3w$. This contradicts the minimality of $[\Lambda :\widetilde{\Lambda}]$, therefore we must have $\Lambda=\widetilde{\Lambda}$. By lemma 3.2.1 below, there is some $x\in K_1$ such that the boundary of the fundamental parallelogram of $q_3\Lambda$ with vertices $x,x+q_3v,x+q_3w,x+q_3(v+w)$ (call it $\mathcal{P}$ again) is a subset of $K_1+q_3\Lambda$. By condition (2), all $K_j$, $j\neq 1$ avoid $K_1+q_3\Lambda$, and hence the boundary of $\mathcal{P}$. Since one translate of $K_j$ by $q_3\Lambda$ intersects $\mathcal{P}$, as it is a fundamental parallelogram of $q_3\Lambda$, \emph{it must lie inside of} $\mathcal{P}$, by convexity (this happens because the boundary of $\mathcal{P}$ splits the plane $\mathbb{R}^2$ into two disjoint regions). Thus, in this case, all $K_j$ for $j>1$ satisfy the additional property
\[\mathfrak{D}K_j\cap q_3(\Lambda\setminus\{0\})=\varnothing.\]
Now let
\[S=\bigl(\bigcup_{j>1}K_j\bigr)\cap \Lambda.\]
From the previous identity we get
\[\mathfrak{D}S\cap q_3(\Lambda\setminus\{0\})=\varnothing,\]
and
\[(K_1-S)\cap q_3\Lambda=\varnothing,\]
from condition (2). Therefore, $K_1$ and $S$ satisfy the conditions of Lemma 2.4, for $r=q_3$, and $d=2$. So, there is a finite set $S'\subset\Lambda$, which is obtained by $S$ by translating each element of $S$ by an element of $q_3\Lambda$, and satisfying
\[\mathfrak{D}S'\cap q_2(\Lambda\setminus\{0\})=\varnothing\]
and
\[(K_1-S')\cap q_2\Lambda=\varnothing,\]
since $q_2\geqslant q_3$. Then,
\begin{multline*}
\sum_{j=1}^n G(K_j,\Lambda)=G(K_1,\Lambda)+\#(S')=\\
=\sum_{r=0}^{q_2-1}\sum_{t\equiv r(\text{mod}q_2)}G(K_1[t]-te^2,\mathbb{Z}e^1)+\sum_{r=0}^{q_2-1}\sum_{t\equiv r(\text{mod}q_2)}\#(S'[t]-te^2).
\end{multline*}
Here, the notation $K[t]$ refers to the original basis $e^1, e^2$. It suffices to prove that for fixed $r$,
\[\sum_{t\equiv r(\text{mod}q_2)}G(K_{1,t},\mathbb{Z}e^1)+\sum_{t\equiv r(\text{mod}q_2)}\#(S'[t]-te^2)\leqslant q_1.\]
We identify $\mathbb{R}e^1$ with $\mathbb{R}$. Hence, we have a finite collection of nonempty compact intervals, $K_{1,t}$, and some lattice points which come from $S'[t]-te^2$. Assume that $S'[t]-te^2=\{m_1e^1,\dotsc,m_ke^1\}$, where $m_1,m_2,\dotsc,m_k$ are distinct integers. Again, we have
\[\mathfrak{D}K_{1,t}\cap q_1(\mathbb{Z}e^1\setminus\{0\})\subset \mathfrak{D}K_j\cap q_1(\Lambda\setminus\{0\})=\varnothing,\]
so condition (1) is satisfied for the intervals $K_{1,t}$ and $m_1e^1,\dotsc,m_ke^1$ (it is trivial for a point). If the intersection
\[(K_{1,t}-K_{1,t'})\cap q_2(\mathbb{Z}e^1)\]
is nonempty for some $t\neq t'$, then there is some $u\in\mathbb{Z}e^1$, such that
\[q_2u+(t-t')e^2\in K_1[t]-K_1[t']\subset \mathfrak{D}K_1,\]
which implies (since $q_2|t-t'$)
\[\mathfrak{D}K_1\cap q_2(\Lambda\setminus\Lambda^1)\neq \varnothing,\]
contradicting condition (1). Furthermore,
\[(K_{1,t}-\{m_ie^1\})\cap q_2(\mathbb{Z}e^1)\subset (K_1-S')\cap q_2\Lambda=\varnothing,\]
and for $i\neq j$,
\[\{m_ie^1\}-\{m_je^1\}\cap q_2(\mathbb{Z}e^1)\subset \mathfrak{D}S'\cap q_2\Lambda=\varnothing,\]
so condition (2) holds as well for the intervals $K_{1,t}$ and the points $m_1e^1,m_2e^1,$ $\dotsc,m_ke^1$, with respect to the lattice $\mathbb{Z}e^1$ and the integers $q_1\geqslant q_2$, hence
\[\sum_{t\equiv r(\text{mod}q_2)}G(K_{1,t},\mathbb{Z}e^1)+\sum_{t\equiv r(\text{mod}q_2)}\#(S'[t]-te^1)\leqslant q_1,\]
as desired, completing the proof.$\square$

%
%
%


\noindent
This implies that inequality (4) is true for $d\leqslant3$. We observe that in order to prove Conjecture 3.2.1 for $d=2$, we used the result for $d=1$. This is exactly the purpose of stating a stronger conjecture than inequality (4); we might be able to use induction on the dimension, something that did not seem possible in this inequality. However, when $d>2$, we need something more than just induction. For $d=2$, Lemma 2.4 was used, because when $D=2$, all but one of the $K_j$ must be confined in a fundamental parallelogram. This is not true in higher dimensions in general; perhaps we need a stronger version of Lemma 2.4.

\noindent
We conclude this section with the following lemma, that was used for the proof of conjecture 3.2.1, case $d=2$:

\noindent
\textbf{\underline{Lemma 3.2.1}} \textit{Let $K\in\mathcal{K}^2$, and $v^1,v^2\in \mathbb{R}^2$ two linearly independent vectors such that the intersections $K\cap(K+v^1)$ and $K\cap(K+v^2)$ are nonempty. Then there exists a point $x\in K$ such that the boundary of the parallelogram with vertices $x,x+v^1,x+v^2,x+v^1+v^2$ is contained in $K+\Lambda$, where $\Lambda$ is the lattice generated by $v^1,v^2$.}

\noindent
\textbf{\underline{Proof}:} From the hypothesis, it is clear that there is a straight line parallel to $v^1$ contained in $K+\mathbb{Z}v^1$, and similarly, a straight line parallel to $v^2$ contained in $K+\mathbb{Z}v^2$. Let $y$ be the point of intersection; then the straight lines parallel to $v^1$, $v^2$, passing through $y$ are contained in $K+\Lambda$. The same happens with any lattice translate of $y$. So, we pick one such translate that belongs to $K$, say $x$. Considering the translates $x+v^1$, $x+v^2$, $x+v^1+v^2$, we deduce that the union of straight lines parallel to $v^1$, $v^2$ and passing through $x$, $x+v^1$, $x+v^2$, $x+v^1+v^2$ is a subset of $K+\Lambda$. It is clear that this union of straight lines contains the boundary of the fundamental parallelogram with vertices $x$, $x+v^1$, $x+v^2$, $x+v^1+v^2$, as desired.$\square$

\section{Reductions of inequality (4)}

\noindent
Two reductions of inequality (4) will be given; the first one is a reduction of conjecture 3.2.1, while the second one is a certain monotonicity property for the discrete measure, that is satisfied by the Lebesgue measure.

\subsection{A simultaneous translation problem}

Observing the proof for the two-dimensional case of conjecture 3.2.1, we see that the main technique was projecting onto a certain hyperplane, and then use induction, i.e., the result for the one-dimensional case. Can we do this in the general case? In particular, what happens when we consider the projections $K_{j,t}=K_j[t]-te^d$ for $1\leqslant j\leqslant n$, $t\equiv r(\text{mod}q_d)$, for a fixed $r$? Do they satisfy conditions (1), (2) of the conjecture, for the lattice $\Lambda^{d-1}$, the basis $e^1,\dotsc,e^{d-1}$ and the integers $q_1\geqslant\dotsb\geqslant q_d$? Not in general. They do, in the special case when $q_{d+1}$ divides $q_d$. If so, we can replace (2) with the weaker condition
\[(K_j-K_l)\cap q_d\Lambda=\varnothing,\]
simply because $q_d\Lambda$ is a sublattice of $q_{d+1}\Lambda$. Indeed,
\[\mathfrak{D}K_{j,t}\cap q_i(\Lambda^{d-1}\setminus\Lambda^{i-1})\subset\mathfrak{D}K_j\cap q_i(\Lambda\setminus\Lambda^{i-1})=\varnothing.\]
For $t\neq t'$, $t\equiv t'(\text{mod}q_d)$, we have
\begin{eqnarray*}
(K_{j,t}-K_{j,t'})\cap q_d\Lambda^{d-1} &=& (K_j[t]-K_j[t'])\cap (q_d\Lambda^{d-1}+(t-t')e^d)\\
&\subset & \mathfrak{D}K_j\cap q_d(\Lambda\setminus\Lambda^{d-1})=\varnothing,
\end{eqnarray*}
and for $j\neq l$, $t\equiv t'(\text{mod}q_d)$, we have
\begin{eqnarray*}
(K_{j,t}-K_{l,t'})\cap q_d\Lambda^{d-1} &=& (K_j[t]-K_l[t'])\cap (q_d\Lambda^{d-1}+(t-t')e^d)\\
&\subset & (K_j-K_l)\cap q_d\Lambda=\varnothing.
\end{eqnarray*}
Hence, as long as $q_{d+1}$ divides $q_d$, we can apply the induction step, using the projection technique. Given the result of conjecture 3.2.1 for $d=2$, we establish the following:

\noindent
\textbf{\underline{Theorem 4.1.1}} \textit{Let $K_1,\dotsc,K_n\subset \mathbb{R}^d$ be convex bodies and $\Lambda\in\mathcal{L}^d$. Also, let $e^1,\dotsc,e^d$ be a basis of $\Lambda$, and denote by $\Lambda^i$ the $\mathbb{Z}$-span of $0,e^1,\dotsc,e^i$, and let $q_1\geqslant q_2\geqslant \dotsb \geqslant q_{d+1}$ be positive integers satisfying
\begin{description}
	\item[(1)] $\mathfrak{D}K_j\cap q_i(\Lambda\setminus\Lambda^{i-1})=\varnothing$ for all $1\leqslant j\leqslant n$ and $1\leqslant i\leqslant d$.
	\item[(2)] $(K_j-K_l)\cap q_{d+1}\Lambda = \varnothing$ for all $1\leqslant j,l\leqslant n$, $j\neq l$.
	\item[(3)] $q_{d+1}|q_d|\dotsb|q_3$.
\end{description}
Then
\[\sum_{j=1}^n G(K_j,\Lambda)\leqslant \prod_{i=1}^d q_i.\]}

\noindent
The purpose is to get rid of the successive divisibility property, (3). What happens when $q_{d+1}$ does not divide $q_d$? We cannot use the same technique anymore, as the projected convex bodies will probably no longer satisfy condition (2). Can we somehow replace $q_{d+1}$ by $q_d$ in condition (2)? We might need to translate the given convex bodies, but we should translate them by a lattice vector, so that the lattice point enumerator remains invariant. We pose the following:

\noindent
\textbf{\underline{Problem}} \textit{Let $K_1,K_2,\dotsc,K_n$ be convex bodies in $\mathbb{R}^d$, $\Lambda$ a lattice, and $r$ be a positive integer, such that the following property holds:
\[(K_i-K_j)\cap r\Lambda=\varnothing,\]
for $i\neq j$, $1\leqslant i,j\leqslant n$. Given a positive integer $t\geqslant r$, is it true that we can translate each $K_i$ by a lattice vector, thus obtaining the convex bodies $K'_1,\dotsc,K'_n$, so that the following property holds for $i\neq j$, $1\leqslant i,j\leqslant n$
\[(K'_i-K'_j)\cap t\Lambda=\varnothing?\]}

\noindent
It is obvious from the analysis at the beginning of the subsection, that if this problem is answered in the affirmative, then it implies conjecture 3.2.1, and consequently inequality (4) for all dimensions. It should be noted that lemma 2.4 is a special case of this problem; also, the case $n=2$ in this problem is a simple consequence of lemma 2.3. Lastly, the one-dimensional case is trivial, or the case where $r$ divides $t$. In this case, we do not have to translate the convex bodies at all.

\noindent
Finally, we state the following corollary to theorem 4.1.1, which is a slight improvement of theorem 3.1.1:

\noindent
\textbf{\underline{Corollary 4.1.1}} \textit{Let $K\in\mathcal{K}^d$, $\Lambda\in\mathcal{L}^d$, $q_i=q_i(K,\Lambda)$. Let $n_1,n_2,\dotsc,n_d$ be a decreasing sequence of positive integers such that
\begin{description}
	\item[(1)] $q_i\leqslant n_i$, for $1\leqslant i\leqslant d$.
	\item[(2)] $n_d|n_{d-1}|\dotsb|n_3$.
\end{description}
Then
\[G(K,\Lambda)\leqslant \prod_{i=1}^{d}n_i.\]
}

\noindent
\textbf{\underline{Proof}:} Let $e^i=e^i(K,\Lambda)$, $\Lambda^i=\Lambda^i(K)$. From the analysis at the beginning of section 3.2, it is clear that the slices $K[t]-te^d$, for $t\equiv r(\text{mod}n_d)$, and numbers $n_1\geqslant n_2\geqslant \dotsb \geqslant n_d$ satisfy conditions (1), (2), and (3) of theorem 4.1.1, whence the desired inequality.$\square$

\noindent
In particular, inequality (4) is verified when $q_d|q_{d-1}|\dotsb|q_3$. This shows that the verification of conjecture 3.2.1 for $d=2$ implies that we need not include the first two terms in this successive divisibility property. And it is clear, that if conjecture 3.2.1 is proven for, say for $d=s$, then inequality (4) is verified when $q_d|q_{d-1}|\dotsb|q_{s+1}$.

\subsection{The discrete monotonicity property}

\noindent
In every proof of Minkowski's second theorem, a monotonicity property for the Lebesgue measure is proven, in one form or another. For example, Bambah [BWZ] proves that
\[\vol(tK/L)\geqslant t^{d-i}\vol(K/L),\]
where $t\geqslant1$, $K\in\mathcal{K}^d$, $L$ a discrete subgroup of $\mathbb{R}^d$ whose rank is equal to $i$, and $\vol(K/L)$ the Lebesgue measure of $K$ taken modulo $L$, i.e., identifying two points of $K$ that are congruent modulo $K$. The above is equivalent to the assertion that
\[\frac{\vol(K/rL)}{r^i}\]
is decreasing in $r>0$. This is the continuous monotonicity property, and holds for all convex bodies $K$ and discrete subgroups $L$ of $\mathbb{R}^d$ unconditionally.

\noindent
Let's now state the discrete monotonicity property; we first replace the Lebesgue measure by a discrete measure, corresponding to a lattice $\Lambda$, so that the measure of a given set $A$ is simply the cardinality of $A\cap\Lambda$. Instead of discrete subgroups of $\mathbb{R}^d$ we consider subgroups of $\Lambda$. Thus:

\noindent
\textbf{\underline{Definition 4.2.1}} \textit{Let $K\in\mathcal{K}^d$, $\Lambda\in\mathcal{L}^d$. We say that $K$ satisfies the \textbf{discrete monotonicity property} with respect to $\Lambda$, if for any subgroup of $\Lambda$, say $\widetilde{\Lambda}$, the following sequence is decreasing in $r>0$, $r\in\mathbb{Z}$:
\[\frac{D_{\Lambda}(K,r\widetilde{\Lambda})}{r^i},\]
where $i$ is the rank of $\widetilde{\Lambda}$.}

\noindent
Here $D_{\Lambda}(K,r\widetilde{\Lambda})$ denotes the cardinality of the set $K\cap\Lambda$ taken modulo $r\widetilde{\Lambda}$. In this case we require that $r$ is an integer, because we always want $r\widetilde{\Lambda}\subset\Lambda$. It is clear that $D_{\Lambda}(K,r\widetilde{\Lambda})$ is the corresponding quantity of $\vol(K/r\Lambda)$ above. Next we prove the following helpful lemma:

\noindent
\textbf{\underline{Lemma 4.2.1}} \textit{Let $K\in\mathcal{K}^d$, $\Lambda\in\mathcal{L}^d$, $a^1,\dotsc,a^d$ $d$ linearly independent vectors of $\Lambda$ and
\[L^i:=\mathbb{Z}a^1\oplus\dotsb\oplus\mathbb{Z}a^i.\]
Assume that $\mathfrak{D}K\cap(L^d\setminus L^i)=\varnothing$. Then
\[D_{\Lambda}(K,L^d)=D_{\Lambda}(K,L^{d-1})=\dotsb=D_{\Lambda}(K,L^i).\]}

\noindent
\textbf{\underline{Proof}:} The hypothesis simply implies that if two points $x,y\in K\cap\Lambda$ are congruent modulo $L^d$, then they must be congruent modulo $L^i$, and consequently congruent modulo $L^j$, for $i\leqslant j\leqslant d$. The lemma then follows from the definition of $D_{\Lambda}(K,L^i)$.$\square$

\noindent
\textbf{\underline{Theorem 4.2.1}} \textit{Assume that $K\in\mathcal{K}^d$ satsifies the discrete monotonicity property with respect to $\Lambda\in\mathcal{L}^d$. Then
\[G(K,\Lambda)\leqslant\prod_{i=1}^{d}q_i(K,\Lambda).\]}

\noindent
\textbf{\underline{Proof}:} Let $\Lambda^i=\Lambda^i(K)$, for $0\leqslant i\leqslant d$, and $q_i=q_i(K,\Lambda)$. By lemma 2.1, we have $\mathfrak{D}K\cap q_i(\Lambda\setminus\Lambda^{i-1})$ for all $i$, and by the virtue of lemma 4.2.1 we have the following series of equalities/inequalities:
\begin{eqnarray*}
q_d^d &\geqslant& D_{\Lambda}(K,q_d\Lambda) = D_{\Lambda}(K,q_d\Lambda^{d-1})\\
&\geqslant& \left(\frac{q_d}{q_{d-1}}\right)^{d-1}D_{\Lambda}(K,q_{d-1}\Lambda^{d-1}) = \left(\frac{q_d}{q_{d-1}}\right)^{d-1}D_{\Lambda}(K,q_{d-1}\Lambda^{d-2})\\
&\vdots&\\
&\geqslant& \left(\frac{q_d}{q_{d-1}}\right)^{d-1}\left(\frac{q_{d-1}}{q_{d-2}}\right)^{d-2} \dotsb \frac{q_2}{q_1}D_{\Lambda}(K,q_1\Lambda^1)\\
&=& \left(\frac{q_d}{q_{d-1}}\right)^{d-1}\left(\frac{q_{d-1}}{q_{d-2}}\right)^{d-2} \dotsb \frac{q_2}{q_1}D_{\Lambda}(K,q_1\Lambda^0)\\
&=& \left(\frac{q_d}{q_{d-1}}\right)^{d-1}\left(\frac{q_{d-1}}{q_{d-2}}\right)^{d-2} \dotsb \frac{q_2}{q_1} G(K,\Lambda)
\end{eqnarray*}
whence
\[G(K,\Lambda)\leqslant \prod_{i=1}^{d}q_i.\square\]

\noindent
The continuous monotonicity property is proven using the homogeneity of the Lebesgue measure. This property is not valid for the discrete measure, so we expect that it might be very difficult to prove the discrete monotonicity property for all convex bodies and all lattices.

\noindent
\textbf{\underline{References}}

\noindent
[BHW] U. Betke, M. Henk, and J.M. Wills, \textit{Successive-minima-type inequalities}, Discrete Comput.
Geom. 9 (1993), no. 2, 165-175.

\noindent
[BWZ] R. P. Bambah, A. C.Woods, and H. Zassenhaus, \textit{Three proofs of Minkowski's second
inequality in the Geometry of Numbers}, J. Austral. Math. Soc. 5 (1965), 453-462.

\noindent
[GL] P.M. Gruber and C.G. Lekkerkerker, \textit{Geometry of Numbers}, 2nd ed., North-Holland, Amsterdam, 1987.

\noindent
[H] M. Henk, \textit{Successive minima and lattice points},  Rendi. Circ. Matematico Palermo, Serie II, Supppl. 70, 2002, 377-384.

\noindent
[M] \underline{$\ \ \ \ \ \ $}, \textit{An optimization problem related to Minkowski's successive minima}, Discrete Comput.
Geom., to appear.

\end{document}